# Statistics of extremes under random censoring

JOHN H.J. EINMAHL[1], AMÉLIE FILS-VILLETARD[2]
and ARMELLE GUILLOU[3]

[1]*Dept. of Econometrics & OR and CentER, Tilburg University, P.O. Box 90153, 5000 LE Tilburg, The Netherlands. E-mail: j.h.j.einmahl@uvt.nl*
[2]*Laboratoire de Statistique Théorique et Appliquée, Université Paris VI, Boîte 158, 4 Place Jussieu, 75252 Paris Cedex 05, France. E-mail: fils@ccr.jussieu.fr*
[3]*IRMA – Département de Mathématiques, Université Louis Pasteur, 7, rue René Descartes, 67084 Strasbourg Cedex, France. E-mail: guillou@math.u-strasbg.fr*

We investigate the estimation of the extreme value index when the data are subject to random censorship. We prove, in a unified way, detailed asymptotic normality results for various estimators of the extreme value index and use these estimators as the main building block for estimators of extreme quantiles. We illustrate the quality of these methods by a small simulation study and apply the estimators to medical data.

## 1. Introduction

Let $X_1, \dots, X_n$ be independent and identically distributed (i.i.d.) random variables, distributed according to an unknown distribution function (df) $F$. A question of great interest is how to obtain a good estimator for a quantile

$$F^{\leftarrow}(1-\varepsilon) = \inf\{y : F(y) \geq 1-\varepsilon\},$$

where $\varepsilon$ is so small that this quantile is situated on the border of, or beyond, the range of the data. Estimating such extreme quantiles is directly linked to the accurate modeling and estimation of the tail of the distribution

$$\overline{F}(x) := 1 - F(x) = \mathbb{P}(X > x)$$

for large thresholds $x$. From extreme value theory, the behaviour of such extreme quantile estimators is known to be governed by one crucial parameter of the underlying distribution, the *extreme value index*. This parameter is important since it measures the tail heaviness of $F$. This estimation has been widely studied in the literature: we mention, for example, Hill (1975), Smith (1987), Dekkers *et al.* (1989) and Drees *et al.* (2004).

However, in classical applications such as the analysis of lifetime data (survival analysis, reliability theory, insurance), a typical feature which appears is censorship. Quite often, $X$ represents the time elapsed from the entry of a patient in, say, a follow-up study until death. If, at the time that the data collection is performed, the patient is still alive or has withdrawn from the study for some reason, the variable of interest $X$ will not be available. A convenient way to model this situation is the introduction of a random variable $Y$, independent of $X$, such that only

$$Z = X \wedge Y \quad \text{and} \quad \delta = \mathbb{1}_{\{X \leq Y\}} \tag{1}$$

are observed. The indicator variable $\delta$ determines whether or not $X$ has been censored. Given a random sample $(Z_i, \delta_i)$, $1 \leq i \leq n$, of independent copies of $(Z, \delta)$, it is our goal to make inference on the right tail of the unknown lifetime distribution function $F$, while $G$, the df of $Y$, is considered to be a nonparametric nuisance parameter.

Statistics of extremes of randomly censored data is a new research field. The statistical problems in this field are difficult since, typically, only a small fraction of the data can be used for inference in the far tail of $F$ and, in the case of censoring, these data are, moreover, not fully informative. The topic was first mentioned in Reiss and Thomas (1997), Section 6.1, where an estimator of a positive extreme value index was introduced, but no (asymptotic) results were derived. Recently, Beirlant *et al.* (2007) proposed estimators for the general extreme value index and for an extreme quantile. That paper made a start on the analysis of the asymptotic properties of some estimators that use the data above a *deterministic* threshold and *only* under the Hall model. In this paper, we consider the "natural" estimators (which are based on the upper order statistics); our methodology is much more general and completely different to their approach.

For almost all applications of extreme value theory, the estimation of the extreme value index is of primary importance. Consequently, it is the main aim of this paper to propose a unified method to prove asymptotic normality for various estimators of the extreme value index under random censoring. We apply our estimators to the problem of extreme quantile estimation under censoring. We illustrate our results with simulations and also apply our methods to AIDS survival data.

We consider data on patients diagnosed with AIDS in Australia before 1 July 1991. The source of these data is Dr P.J. Solomon and the Australian National Centre in HIV Epidemiology and Clinical Research; see Venables and Ripley (2002). The information on each patient includes gender, date of diagnosis, date of death or end of observation and an indicator as to which of the two is the case. The data set contains 2843 patients, of which 1761 died; the other survival times are right-censored. We will apply our methodology to the 2754 male patients (there are only 89 women in the data set), of which 1708 died. Apart from assessing the heaviness of the right tail of the survival function $1 - F$ by means of the estimation of the extreme value index, it is also important to estimate very high quantiles of $F$, thus obtaining a good indication of how long very strong men will survive AIDS.

Another possible application, not pursued in this paper, is to annuity insurance contracts. Life annuities are contractual guarantees, issued by insurance companies, pension plans and government retirement systems, that offer promises to provide periodic income over the lifetime of individuals. If we monitor the policyholders during a certain period, the data are right censored since many policyholders survive until the end of the observation period. We are interested in the far right tail of the future lifetime distribution of the annuitants, since longevity is an important and difficult risk to evaluate for insurance companies. In the case of life annuities, it needs to be estimated as accurately as possible for setting adequate insurance premiums.

We will study estimators for the extreme value index of $F$, assuming that $F$ and $G$ are both in the max-domain of attraction of an extreme value distribution. In Section 2, we introduce various estimators of this extreme value index and we establish, in a unified way, their asymptotic behaviors. We also introduce an estimator for very high quantiles. Various examples are given in Section 3 and a small simulation study is performed. Our estimators are applied to the AIDS data in Section 4.

## 2. Estimators and main results

Let $X_1, \dots, X_n$ be a sequence of i.i.d. random variables from a df $F$. We denote the order statistics by

$$X_{1,n} \leq \cdots \leq X_{n,n}.$$

The weak convergence of the centered and standardized maxima $X_{n,n}$ implies the existence of sequences of constants $a_n > 0$ and $b_n$ and a df $\widetilde{G}$ such that

$$\lim_{n \to \infty} \mathbb{P}\left( \frac{X_{n,n} - b_n}{a_n} \leq x \right) = \widetilde{G}(x) \qquad (2)$$

for all $x$ where $\widetilde{G}$ is continuous. The work of Fisher and Tippett (1928), Gnedenko (1943) and de Haan (1970) answered the question on the possible limits and characterized the classes of distribution functions $F$ having a certain limit in (2).

This convergence result is our main assumption. Up to location and scale, the possible limiting dfs $\widetilde{G}$ in (2) are given by the so-called *extreme value distributions* $G_\gamma$, defined by

$$G_\gamma(x) = \begin{cases} \exp(-(1 + \gamma x)^{-1/\gamma}), & \text{if } \gamma \neq 0, \\ \exp(-\exp(-x)), & \text{if } \gamma = 0. \end{cases} \qquad (3)$$

We say that $F$ is in the (max-)domain of attraction of $G_\gamma$, denoting this by $F \in D(G_\gamma)$. Here $\gamma$ is the extreme value index. Knowledge of $\gamma$ is crucial for estimating the right tail of $F$.

We briefly review some estimators of $\gamma$ that have been proposed in the literature. The most famous is probably the Hill (1975) estimator

$$\widehat{\gamma}_{X,k,n}^{(H)} := M_{X,k,n}^{(1)} := \frac{1}{k} \sum_{i=1}^{k} \log X_{n-i+1,n} - \log X_{n-k,n}, \qquad (4)$$

where $k \in \{1, \ldots, n-1\}$. However, this estimator is only useful when $\gamma > 0$. A generalization which works for any $\gamma \in \mathbb{R}$ is the so-called *moment estimator*, introduced in Dekkers *et al.* (1989):

$$\widehat{\gamma}_{X,k,n}^{(M)} := M_{X,k,n}^{(1)} + S_{X,k,n} := M_{X,k,n}^{(1)} + 1 - \frac{1}{2}\left(1 - \frac{(M_{X,k,n}^{(1)})^2}{M_{X,k,n}^{(2)}}\right)^{-1}, \tag{5}$$

with

$$M_{X,k,n}^{(2)} := \frac{1}{k}\sum_{i=1}^{k}(\log X_{n-i+1,n} - \log X_{n-k,n})^2.$$

The Hill estimator can be derived in several ways, a very appealing one being the slope of the Pareto quantile plot, which consists of the points

$$\left(\log \frac{n+1}{i}, \log X_{n-i+1,n}\right), \qquad i = 1, \ldots, k.$$

This plot has been generalized in Beirlant *et al.* (1996) by defining $UH_{i,n} = X_{n-i,n}\widehat{\gamma}_{X,i,n}^{(H)}$ and considering the points

$$\left(\log \frac{n+1}{i}, \log UH_{i,n}\right), \qquad i = 1, \ldots, k.$$

This generalized quantile plot becomes almost linear for small enough $k$, that is, for extreme values. It follows immediately that the slope of this graph will estimate $\gamma$ regardless of whether it is positive, negative or zero. An estimator of this slope is given by

$$\widehat{\gamma}_{X,k,n}^{(UH)} := \frac{1}{k}\sum_{i=1}^{k}\log UH_{i,n} - \log UH_{k+1,n}, \tag{6}$$

where $k \in \{1, \ldots, n-2\}$.

A quite different estimator of $\gamma$ is the so-called *maximum likelihood* (*ML*) estimator. (Note that the classical, parametric ML approach is not applicable because $F$ is not in a parametric family.) The approach relies on results in Balkema and de Haan (1974) and Pickands (1975), stating that the limit distribution of the exceedances $E_j = X_j - t$ ($X_j > t$) over a threshold $t$, when $t$ tends to the right end-point of $F$, is given by a generalized Pareto distribution depending on two parameters, $\gamma$ and $\sigma$. In practice, $t$ is replaced by an order statistic $X_{n-k,n}$ and the resulting *ML*-estimators are denoted by $\widehat{\gamma}_{X,k,n}^{(ML)}$ and $\widehat{\sigma}_{X,k,n}^{(ML)}$.

In the case of censoring, we would like to adapt all of these methods. Actually, we will provide a general adaptation of estimators of the extreme value index and a unified proof of their asymptotic normality; the four estimators above are special cases of this. We assume that both $F$ and $G$ are absolutely continuous and that $F \in D(G_{\gamma_1})$ and $G \in D(G_{\gamma_2})$ for some $\gamma_1, \gamma_2 \in \mathbb{R}$. The extreme value index of $H$, the df of $Z$ defined in (1), exists and is denoted by $\gamma$. Let $\tau_F = \sup\{x : F(x) < 1\}$ (resp., $\tau_G$ and $\tau_H$) denote the right endpoint of the support of $F$ (resp., $G$ and $H$). In the sequel, we assume that the pair $(F, G)$ is in one of the following three cases:

$$\begin{cases} \text{case 1: } \gamma_1 > 0, \gamma_2 > 0, & \text{in this case } \gamma = \dfrac{\gamma_1 \gamma_2}{\gamma_1 + \gamma_2}, \\ \text{case 2: } \gamma_1 < 0, \gamma_2 < 0, \tau_F = \tau_G, & \text{in this case } \gamma = \dfrac{\gamma_1 \gamma_2}{\gamma_1 + \gamma_2}, \\ \text{case 3: } \gamma_1 = \gamma_2 = 0, \tau_F = \tau_G = \infty, & \text{in this case } \gamma = 0. \end{cases} \quad (7)$$

(In case 3, we also define, for convenient presentation, $\frac{\gamma_1 \gamma_2}{\gamma_1 + \gamma_2} = \gamma = 0$.) The other possibilities are not very interesting. Typically, they are very close to the "uncensored case", which has been studied in detail in the literature (this holds, in particular, when $\gamma_1 < 0$ and $\gamma_2 > 0$) or the "completely censored situation", where estimation is impossible (this holds, in particular, when $\gamma_1 > 0$ and $\gamma_2 < 0$).

The first important point that should be mentioned is the fact that all of the preceding estimators (Hill, moment, *UH* or *ML*) are obviously not consistent if they are based on the sample $Z_1, \dots, Z_n$, that is, if the censoring is not taken into account. Indeed, they all converge to $\gamma$, the extreme value index of the $Z$-sample, and not to $\gamma_1$, the extreme value index of $F$. Consequently, we must adapt all of these estimators to censoring. We will divide all these estimators by the proportion of non-censored observations in the $k$ largest $Z$'s:

$$\widehat{\gamma}_{Z,k,n}^{(c,\cdot)} = \frac{\widehat{\gamma}_{Z,k,n}^{(\cdot)}}{\widehat{p}}, \qquad \text{where } \widehat{p} = \frac{1}{k} \sum_{j=1}^{k} \delta_{[n-j+1,n]},$$

with $\delta_{[1,n]}, \dots, \delta_{[n,n]}$ being the $\delta$'s corresponding to $Z_{1,n}, \dots, Z_{n,n}$, respectively. $\widehat{\gamma}_{Z,k,n}^{(\cdot)}$ could be any estimator not adapted to censoring, in particular, $\widehat{\gamma}_{Z,k,n}^{(H)}, \widehat{\gamma}_{Z,k,n}^{(M)}, \widehat{\gamma}_{Z,k,n}^{(UH)}$ or $\widehat{\gamma}_{Z,k,n}^{(ML)}$. It will follow that $\widehat{p}$ estimates $\frac{\gamma_2}{\gamma_1 + \gamma_2}$, hence $\widehat{\gamma}_{Z,k,n}^{(\cdot)}$ estimates $\gamma$ divided by $\frac{\gamma_2}{\gamma_1 + \gamma_2}$, which is equal to $\gamma_1$. It is our main aim to study in detail the asymptotic normality of these estimators.

To illustrate the difference between the estimators, adapted and not adapted to censoring, in Figure 1(a), we plot $\widehat{\gamma}_{Z,k,n}^{(UH)}$ (dashed line) and $\widehat{\gamma}_{Z,k,n}^{(c,\cdot)}$ (full line) as a function of $k$ for the AIDS survival data. We see a quite stable plot when $k$ ranges from about 200 (or 350) to 1200 and a substantial difference between the two estimators. Similar graphs could be presented for the other estimators.

Let us now consider the estimation of an extreme quantile $x_\varepsilon = F^{\leftarrow}(1 - \varepsilon)$. Denoting by $\widehat{F}_n$ the Kaplan–Meier (1958) product-limit estimator, we can adapt the classical estimators proposed in the literature as follows:

$$\widehat{x}_{\varepsilon,k}^{(c,\cdot)} = Z_{n-k,n} + \widehat{a}_{Z,k,n}^{(c,\cdot)} \frac{\left( (1 - \widehat{F}_n(Z_{n-k,n})) / \varepsilon \right)^{\widehat{\gamma}_{Z,k,n}^{(c,\cdot)}} - 1}{\widehat{\gamma}_{Z,k,n}^{(c,\cdot)}}, \quad (8)$$

where

$$\widehat{a}_{Z,k,n}^{(c,\cdot)} = \frac{Z_{n-k,n} M_{Z,k,n}^{(1)}(1 - S_{Z,k,n})}{\widehat{p}} \qquad \text{for } M \text{ and } UH,$$

with $S_{Z,k,n}$ defined in (5) and

$$\widehat{a}_{Z,k,n}^{(c,ML)} = \frac{\widehat{\sigma}_{Z,k,n}^{(ML)}}{\widehat{p}}.$$

Note that these estimators are defined under the assumption that the two endpoints $\tau_F$ and $\tau_G$ are equal, but possibly infinite. This is true for the three cases defined in (7). Also, note that we have excluded the Hill estimator since it only works in case 1.

Again, to observe the difference between the adapted and non-adapted estimators, in Figure 1(b), we plot $\widehat{x}_{0.001,k}^{(UH)}$ (dashed line) and $\widehat{x}_{0.001,k}^{(c,UH)}$ (full line) for the AIDS data. The difference between the two estimators (for $k$ between 250 and 500) is about 10 years.

Beirlant *et al.* (2007) considered asymptotic properties of *some* of these estimators, when $Z_{n-k,n}$ is replaced by a deterministic $t$ in the preceding formulas and only under the Hall model. Also, note that the asymptotic bias of these estimators has not been studied. Our aim in this paper is to establish the asymptotic normality (including bias and variance) of all of the above estimators of the extreme value index (based on $k+1$ upper order statistics or, equivalently, on a random threshold $Z_{n-k,n}$). We use a general approach that separates extreme value theory and censoring. Therefore, in the proof, we can treat the above four estimators (and others) simultaneously.

To specify the asymptotic bias of the different estimators, we use a second-order condition phrased in terms of the tail quantile function $U_H(x) = H^{\leftarrow}(1 - \frac{1}{x})$. From the theory of generalized regular variation of second-order outlined in de Haan and Stadtmüller (1996), we assume the existence of a positive function $a$ and a second eventually positive function $a_2$ with $\lim_{x \to \infty} a_2(x) = 0$, such that the limit

$$\lim_{x \to \infty} \frac{1}{a_2(x)} \left\{ \frac{U_H(ux) - U_H(x)}{a(x)} - h_\gamma(u) \right\} = k(u) \tag{9}$$

exists for $u \in (0, \infty)$, with $h_\gamma(u) = \int_1^u z^{\gamma-1} \, \mathrm{d}z$. It follows that there exists a $c \in \mathbb{R}$ and a second-order parameter $\rho \le 0$ for which the function $a$ satisfies

$$\lim_{x \to \infty} \left\{ \frac{a(ux)}{a(x)} - u^\gamma \right\} \Big/ a_2(x) = c u^\gamma h_\rho(u). \tag{10}$$

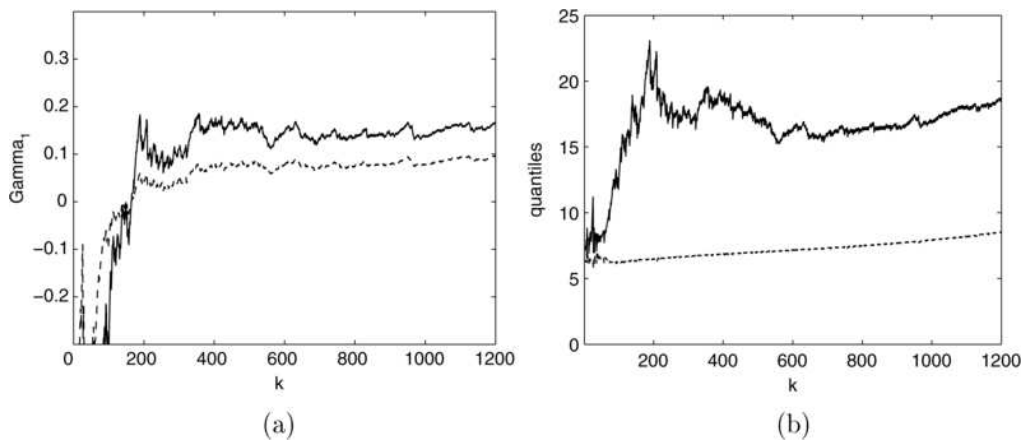

(a)

(b)

**Figure 1.** *UH*-estimator adapted (full line) and not adapted (dashed line) to censoring (a) for the extreme value index and (b) for the extreme quantile with $\varepsilon = 0.001$ for the AIDS survival data.

The function $a_2$ is regularly varying with index $\rho$. As usual, we will assume that $\rho < 0$ and we will also assume that the slowly varying part of $a_2$ is asymptotically equivalent to a positive constant, which can and will always be taken equal to 1. For an appropriate choice of the function $a$, the function $k$ that appears in (9) admits the representation

$$k(u) = A h_{\gamma+\rho}(u), \tag{11}$$

with $A \ne 0$; $c$ in (10) is now equal to 0. We denote the class of second-order regularly varying functions $U_H$ (satisfying (9)–(11) with $c = 0$) by $GRV_2(\gamma, \rho; a(x), a_2(x); A)$.

From Vanroelen (2003), we obtain the following representations of $U_H$ (see also the Appendix in Draisma *et al.* (1999)):

- $0 < -\rho < \gamma$: for $U_H \in GRV_2(\gamma, \rho; \ell_+ x^\gamma, a_2(x); A)$,

$$U_H(x) = \ell_+ x^\gamma \left\{ \frac{1}{\gamma} + \frac{A}{\gamma + \rho} a_2(x)(1 + o(1)) \right\};$$

- $\gamma = -\rho$: for $U_H \in GRV_2(\gamma, -\gamma; \ell_+ x^\gamma, x^{-\gamma} \ell_2(x); A)$,

$$U_H(x) = \ell_+ x^\gamma \left\{ \frac{1}{\gamma} + x^{-\gamma} L_2(x) \right\},$$

with $L_2(x) = B + \int_1^x (A + o(1)) \frac{\ell_2(t)}{t} \, \mathrm{d}t + o(\ell_2(x))$ for some constant $B$ and some slowly varying function $\ell_2$;

- $0 < \gamma < -\rho$: for $U_H \in GRV_2(\gamma, \rho; \ell_+ x^\gamma, a_2(x); A)$,

$$U_H(x) = \ell_+ x^\gamma \left\{ \frac{1}{\gamma} + D x^{-\gamma} + \frac{A}{\gamma + \rho} a_2(x)(1 + o(1)) \right\}$$

  (so $D = \frac{1}{\ell_+} \lim_{x \to \infty} \{ U_H(x) - a(x)/\gamma \}$);

- $\gamma = 0$: for $U_H \in GRV_2(0, \rho; \ell_+, a_2(x); A)$,

$$U_H(x) = \ell_+ \log x + D + \frac{A \ell_+}{\rho} a_2(x)(1 + o(1));$$

- $\gamma < 0$: for $U_H \in GRV_2(\gamma, \rho; \ell_+ x^\gamma, a_2(x); A)$,

$$U_H(x) = \tau_H - \ell_+ x^\gamma \left\{ \frac{1}{-\gamma} - \frac{A}{\gamma + \rho} a_2(x)(1 + o(1)) \right\},$$

where $\ell_+ > 0, A \neq 0, D \in \mathbb{R}$.

In the statement of our results, we use the following notation, similar to that used in Beirlant *et al.* (2005):

$$b(x) = \begin{cases} \dfrac{A\rho[\rho + \gamma(1 - \rho)]}{(\gamma + \rho)(1 - \rho)} a_2(x), & \text{if } 0 < -\rho < \gamma \text{ or if } 0 < \gamma < -\rho \text{ with } D = 0, \\[2mm] -\dfrac{\gamma^3}{(1 + \gamma)} x^{-\gamma} L_2(x), & \text{if } \gamma = -\rho, \\[2mm] -\dfrac{\gamma^3 D}{(1 + \gamma)} x^{-\gamma}, & \text{if } 0 < \gamma < -\rho \text{ with } D \neq 0, \\[2mm] \dfrac{1}{\log^2 x}, & \text{if } \gamma = 0, \\[2mm] \dfrac{A\rho(1 - \gamma)}{(1 - \gamma - \rho)} a_2(x), & \text{if } \gamma < \rho, \\[2mm] -\dfrac{\gamma}{1 - 2\gamma} \dfrac{\ell_+}{\tau_H} x^\gamma, & \text{if } \rho < \gamma < 0, \\[2mm] \dfrac{\gamma}{1 - 2\gamma} \left[ A(1 - \gamma) - \dfrac{\ell_+}{\tau_H} \right] x^\gamma, & \text{if } \gamma = \rho \end{cases}$$

and

$$\widetilde{\rho} = \begin{cases} -\gamma, & \text{if } 0 < \gamma < -\rho \text{ with } D \neq 0, \\ \rho, & \text{if } -\rho \leq \gamma \text{ or if } 0 < \gamma < -\rho \text{ with } D = 0, \text{ or if } \gamma < \rho, \\ \gamma, & \text{if } \rho \leq \gamma \leq 0. \end{cases}$$

Before stating our main result, define

$$p(z) = \mathbb{P}(\delta = 1 | Z = z).$$

It follows that

$$p(z) = \frac{(1 - G(z))f(z)}{(1 - G(z))f(z) + (1 - F(z))g(z)}, \tag{12}$$

where $f$ and $g$ denote the densities of $F$ and $G$, respectively. Note that, in cases 1 and 2, $\lim_{z \to \tau_H} p(z)$ exists and is equal to $\frac{\gamma_2}{\gamma_1 + \gamma_2} =: p \in (0, 1)$. Assume that, in case 3, this limit also exists and is positive and again denote it by $p$. By convention, we also define $\frac{\gamma_2}{\gamma_1 + \gamma_2} = p$ for that case.

In the sequel, $k = k_n$ is an intermediate sequence, that is, a sequence such that $k \to \infty$ and $\frac{k}{n} \to 0$, as $n \to \infty$. Our main result now reads as follows.

**Theorem 1.** *Under the assumptions that, for $n \to \infty$,*

$$\begin{cases} \sqrt{k}\,a_2\left(\dfrac{n}{k}\right) \to \alpha_1 \in \mathbb{R}, & \text{for the ML-estimator,} \\[2mm] \sqrt{k}\,b\left(\dfrac{n}{k}\right) \to \alpha_1 \in \mathbb{R}, & \text{for the other three estimators,} \end{cases} \tag{13}$$

$$\frac{1}{\sqrt{k}} \sum_{i=1}^{k} \left[ p\left(H^{\leftarrow}\left(1 - \frac{i}{n}\right)\right) - p \right] \longrightarrow \alpha_2 \in \mathbb{R} \tag{14}$$

*and*

$$\sqrt{k} \sup_{\{1-k/n \leq t < 1, |t-s| \leq C\sqrt{k}/n, s < 1\}} |p(H^{\leftarrow}(t)) - p(H^{\leftarrow}(s))| \longrightarrow 0 \quad \text{for all } C > 0, \tag{15}$$

*we have, for the four estimators (for the Hill estimator, we assume case 1 holds and for the ML-estimator, that $\gamma > -\frac{1}{2}$),*

$$\sqrt{k}(\widehat{\gamma}_{Z,k,n}^{(c,\cdot)} - \gamma_1) \xrightarrow{d} \mathcal{N}\left( \frac{1}{p}(\alpha_1 b_0 - \gamma_1 \alpha_2), \frac{\sigma^2 + \gamma_1^2 p(1-p)}{p^2} \right),$$

*where $\alpha_1 b_0$ (resp., $\sigma^2$) denotes the bias (resp., the variance) of $\sqrt{k}(\widehat{\gamma}_{Z,k,n}^{(\cdot)} - \gamma)$.*

This leads to the following corollary, the proof of which is rather straightforward. For the Hill estimator, the asymptotic bias-term follows easily from direct computations and for the other three estimators, it follows from the expressions for the asymptotic bias-terms of the corresponding "uncensored" estimators: see Beirlant *et al.* (2005) and Drees *et al.* (2004).

**Corollary 1.** *Under the assumptions of Theorem 1, we have*

$$\sqrt{k}(\widehat{\gamma}_{Z,k,n}^{(c,H)} - \gamma_1)$$

$$\xrightarrow{d} \mathcal{N}\left( \mu^{(c,H)}, \frac{\gamma_1^3}{\gamma} \right) \quad \text{in case 1;}$$

$$\sqrt{k}(\widehat{\gamma}_{Z,k,n}^{(c,M)} - \gamma_1)$$

$$\xrightarrow{d} \begin{cases} \mathcal{N}\left( \mu^{(c,M)}, \dfrac{\gamma_1^2}{\gamma^2}(1 + \gamma_1 \gamma) \right), & \text{in case 1,} \\[3mm] \mathcal{N}\left( \mu^{(c,M)}, \dfrac{\gamma_1^2(1-\gamma)^2(1-2\gamma)(1-\gamma+6\gamma^2)}{\gamma^2(1-4\gamma)(1-3\gamma)} + \gamma_1^2\left(\dfrac{\gamma_1}{\gamma} - 1\right) \right), & \text{in case 2,} \\[3mm] \mathcal{N}(\mu^{(c,M)}, p^{-2}), & \text{in case 3;} \end{cases}$$

$$\sqrt{k}(\widehat{\gamma}_{Z,k,n}^{(c,UH)} - \gamma_1)$$

$$\overset{d}{\longrightarrow} \begin{cases} \mathcal{N}\left(\mu^{(c,UH)}, \dfrac{\gamma_1^2}{\gamma^2}(1+\gamma_1\gamma)\right), & \text{in case 1,} \\[2ex] \mathcal{N}\left(\mu^{(c,UH)}, \dfrac{\gamma_1^2(1-\gamma)(1+\gamma+2\gamma^2)}{\gamma^2(1-2\gamma)} + \gamma_1^2\left(\dfrac{\gamma_1}{\gamma}-1\right)\right), & \text{in case 2,} \\[2ex] \mathcal{N}(\mu^{(c,UH)}, p^{-2}), & \text{in case 3;} \end{cases}$$

$$\sqrt{k}(\widehat{\gamma}_{Z,k,n}^{(c,ML)} - \gamma_1)$$

$$\overset{d}{\longrightarrow} \mathcal{N}(\mu^{(c,ML)}, p^{-2}[1+\gamma(2+\gamma_1)]), \qquad \text{in cases 1, 3 and 2 with } \gamma > -\tfrac{1}{2},$$

*where*

$$\mu^{(c,H)} := -\frac{\gamma_1\alpha_2}{p} + \frac{\alpha_1}{p}\frac{\gamma}{\widehat{\rho}+\gamma(1-\widehat{\rho})};$$

$$\mu^{(c,M)} := -\frac{\gamma_1\alpha_2}{p}$$

$$+\frac{\alpha_1}{p} \cdot \begin{cases} \dfrac{1}{1-\widehat{\rho}}, & \text{in case 1,} \\[2ex] \dfrac{2\gamma-1}{\widehat{\rho}(1-\widehat{\rho})}, & \text{in case 2, if } \rho < \gamma, \\[2ex] \dfrac{1-2\gamma}{(1-\gamma)(1-3\gamma)} \dfrac{A(1-\gamma)^2 - (\gamma+1)\frac{\ell_+}{\tau_H}}{A(1-\gamma) - \frac{\ell_+}{\tau_H}}, & \text{in case 2, if } \rho = \gamma, \\[2ex] \dfrac{1-2\gamma}{1-2\gamma-\widehat{\rho}}, & \text{in case 2, if } \gamma < \rho, \\[1ex] 1, & \text{in case 3;} \end{cases}$$

$$\mu^{(c,UH)} := -\frac{\gamma_1\alpha_2}{p} + \frac{\alpha_1}{p(1-\widehat{\rho})};$$

$$\mu^{(c,ML)} := -\frac{\gamma_1\alpha_2}{p} + \frac{\alpha_1}{p}\frac{\rho(\gamma+1)A}{(1-\rho)(1-\rho+\gamma)}.$$

**Proof of Theorem 1.** We consider the following decomposition

$$\sqrt{k}(\widehat{\gamma}_{Z,k,n}^{(c,\cdot)} - \gamma_1) = \frac{1}{\widehat{p}}\sqrt{k}(\widehat{\gamma}_{Z,k,n}^{(\cdot)} - \gamma) + \frac{1}{\widehat{p}}\sqrt{k}(\gamma - \gamma_1\widehat{p})$$

$$= \frac{1}{\widehat{p}}\sqrt{k}(\widehat{\gamma}_{Z,k,n}^{(\cdot)} - \gamma) + \frac{\gamma_1}{\widehat{p}}\sqrt{k}\left(\frac{\gamma_2}{\gamma_1+\gamma_2} - \widehat{p}\right). \tag{16}$$

The asymptotic behavior of $\sqrt{k}(\widehat{\gamma}_{Z,k,n}^{(\cdot)} - \gamma)$ is well known since this estimator is based on the $Z$-sample, that is, on the uncensored situation; see Beirlant *et al.* (2005) and Drees *et al.* (2004).

First, note that in case 3, $\gamma_1 = \gamma = 0$. Therefore, the second term in the decomposition (16) is exactly 0 provided $\widehat{p} > 0$. That means that this case follows, since $\widehat{p} \overset{\mathbb{P}}{\longrightarrow} p > 0$. We now focus in detail on the second term of the decomposition in (16) for the cases 1 and 2.

To this end, consider the following construction. Let $Z$ be a random variable with df $H$. Let $U$ have a uniform$(0,1)$ distribution and be independent of $Z$. Define

$$\delta = \begin{cases} 1, & \text{if } U \leq p(Z), \\ 0, & \text{if } U > p(Z) \end{cases}$$

and

$$\widetilde{\delta} = \begin{cases} 1, & \text{if } U \leq p, \\ 0, & \text{if } U > p. \end{cases}$$

We repeat this construction independently $n$ times. It is easy to show that the resulting pairs $(Z_i, \delta_i)$, $i = 1, \ldots, n$, have the same distribution as the initial pairs $(Z_i, \delta_i)$, $i = 1, \ldots, n$, for all $n \in \mathbb{N}$, so we continue with the new pairs $(Z_i, \delta_i)$.

Moreover, $Z$ and $\widetilde{\delta}$ are clearly independent and satisfy

$$\mathbb{P}(|\delta - \widetilde{\delta}| = 1 | Z = z) = |p - p(z)|.$$

Consider the order statistics $Z_{1,n} \leq \cdots \leq Z_{n,n}$ and denote the induced order statistics of the $U$'s by $U_{[1,n]}, \ldots, U_{[n,n]}$. We can write $\widehat{p}$ as follows:

$$\widehat{p} = \frac{1}{k} \sum_{j=1}^{k} \mathbb{1}_{\{U_{[n-j+1,n]} \leq p(Z_{n-j+1,n})\}}$$

and, similarly,

$$\widetilde{p} := \frac{1}{k} \sum_{j=1}^{k} \widetilde{\delta}_{[n-j+1,n]} = \frac{1}{k} \sum_{j=1}^{k} \mathbb{1}_{\{U_{[n-j+1,n]} \leq p\}}.$$

Clearly, $U_{[1,n]}, \ldots, U_{[n,n]}$ are i.i.d. and independent of the $Z$-sample.

We use the following decomposition:

$$\sqrt{k}(\widehat{p} - p) = \sqrt{k}(\widehat{p} - \widetilde{p}) + \sqrt{k}(\widetilde{p} - p). \tag{17}$$

Since $\widetilde{p} \overset{d}{=} \frac{1}{k} \sum_{j=1}^{k} \mathbb{1}_{\{U_j \leq p\}}$, we have

$$\sqrt{k}(\widetilde{p} - p) \overset{d}{\longrightarrow} \mathcal{N}(0, p(1-p)).$$

Now, we are interested in $\sqrt{k}(\widehat{p} - \widetilde{p})$, which turns out to be a bias term. It can be rewritten as follows:

$$\sqrt{k}(\widehat{p} - \widetilde{p}) \overset{d}{=} \frac{1}{\sqrt{k}} \sum_{j=1}^{k} [\mathbb{1}_{\{U_j \leq p(Z_{n-j+1,n})\}} - \mathbb{1}_{\{U_j \leq p\}}]$$

$$= \frac{1}{\sqrt{k}} \sum_{j=1}^{k} [\mathbb{1}_{\{U_j \le p(Z_{n-j+1,n})\}} - \mathbb{1}_{\{U_j \le p(H^{\leftarrow}(1-j/n))\}}]$$

$$+ \frac{1}{\sqrt{k}} \sum_{j=1}^{k} [\mathbb{1}_{\{U_j \le p(H^{\leftarrow}(1-j/n))\}} - \mathbb{1}_{\{U_j \le p\}}]$$

$$=: T_{1,k} + T_{2,k}.$$

Under the assumptions (14) and (15), the convergence in probability of $T_{2,k}$ to $\alpha_2$ then follows from a result in Chow and Teicher (1997), page 356.

So we now need to show that $T_{1,k} \xrightarrow{\mathbb{P}} 0$. To this end, write $V_i = H(Z_i)$ so that $Z_i = H^{\leftarrow}(V_i)$. The $V_i$ are i.i.d. uniform $(0,1)$. Also, write $r(t) = p(H^{\leftarrow}(t))$. Then

$$T_{1,k} = \frac{1}{\sqrt{k}} \sum_{j=1}^{k} [\mathbb{1}_{\{U_j \le r(V_{n-j+1,n})\}} - \mathbb{1}_{\{U_j \le r(1-j/n)\}}].$$

By the weak convergence of the uniform tail quantile process, we have, uniformly in $1 \le j \le k$,

$$V_{n-j+1,n} - \left(1 - \frac{j}{n}\right) = O_{\mathbb{P}}\left(\frac{\sqrt{k}}{n}\right).$$

Let $\eta > 0$. Using (15), we have, with arbitrarily high probability, for large $n$,

$$|T_{1,k}| \le \frac{1}{\sqrt{k}} \sum_{j=1}^{k} |\mathbb{1}_{\{U_j \le r(V_{n-j+1,n})\}} - \mathbb{1}_{\{U_j \le r(1-j/n)\}}|$$

$$\stackrel{d}{=} \frac{1}{\sqrt{k}} \sum_{j=1}^{k} \mathbb{1}_{\{U_j \le |r(V_{n-j+1,n}) - r(1-j/n)|\}}$$

$$\le \frac{1}{\sqrt{k}} \sum_{j=1}^{k} \mathbb{1}_{\{U_j \le \eta/\sqrt{k}\}}.$$

Using the aforementioned result in Chow and Teicher (1997), page 356, and the fact that $\eta > 0$ can be chosen arbitrarily small, $T_{1,k} \xrightarrow{\mathbb{P}} 0$ follows.

Finally, combining (16) and (17) yields

$$\sqrt{k}(\widehat{\gamma}_{Z,k,n}^{(c,\cdot)} - \gamma_1) = \frac{1}{\widehat{p}}(\sqrt{k}(\widehat{\gamma}_{Z,k,n}^{(\cdot)} - \gamma) - \gamma_1\sqrt{k}(\widetilde{p} - p)) - \frac{\gamma_1\alpha_2}{\widehat{p}} + o_{\mathbb{P}}(1), \qquad (18)$$

with the two terms within the brackets independent since the first is based on the $Z$-sample and the second on the $U$-sample. Therefore, under the assumptions (13)–(15), we

have

$$\sqrt{k}(\hat{\gamma}_{Z,k,n}^{(c,\bullet)} - \gamma_1) \xrightarrow{d} \mathcal{N}\left(\frac{1}{p}(\alpha_1 b_0 - \gamma_1 \alpha_2), \frac{\sigma^2 + \gamma_1^2 p(1-p)}{p^2}\right). \qquad \square$$

## 3. Examples and small simulation study

In this section, we consider three examples: first, a Burr distribution censored by another Burr distribution (hence an example of case 1), second a reverse Burr distribution censored by another reverse Burr distribution (an example of case 2) and finally a logistic distribution censored by a logistic distribution (case 3). We show that these distributions satisfy all of the assumptions and calculate the bias terms explicitly. In particular, we will see how assumptions (13) and (14) compare. We also provide simulations to illustrate the behavior of our estimators for these distributions.

**Example 1.** $X \sim \mathrm{Burr}(\beta_1, \tau_1, \lambda_1)$ and $Y \sim \mathrm{Burr}(\beta_2, \tau_2, \lambda_2)$, $\beta_1, \tau_1, \lambda_1, \beta_2, \tau_2, \lambda_2 > 0$.

In that case,

$$1 - F(x) = \left(\frac{\beta_1}{\beta_1 + x^{\tau_1}}\right)^{\lambda_1}$$

$$= x^{-\tau_1 \lambda_1} \beta_1^{\lambda_1} (1 + \beta_1 x^{-\tau_1})^{-\lambda_1}, \qquad x > 0;$$

$$1 - G(x) = \left(\frac{\beta_2}{\beta_2 + x^{\tau_2}}\right)^{\lambda_2}$$

$$= x^{-\tau_2 \lambda_2} \beta_2^{\lambda_2} (1 + \beta_2 x^{-\tau_2})^{-\lambda_2}, \qquad x > 0.$$

We can infer that

$$U_H(x) = H^{\leftarrow}\left(1 - \frac{1}{x}\right)$$

$$= (\beta_1^{\lambda_1} \beta_2^{\lambda_2} x)^{1/(\tau_1 \lambda_1 + \tau_2 \lambda_2)} [1 - \gamma \eta (\beta_1^{\lambda_1} \beta_2^{\lambda_2} x)^{\rho} (1 + o(1))],$$

with

$$\tau = \min(\tau_1, \tau_2), \qquad \rho = -\gamma \tau$$

and

$$\eta = \begin{cases} \lambda_1 \beta_1, & \text{if } \tau_1 < \tau_2, \\ \lambda_2 \beta_2, & \text{if } \tau_1 > \tau_2, \\ \lambda_1 \beta_1 + \lambda_2 \beta_2, & \text{if } \tau_1 = \tau_2. \end{cases}$$

The parameters of interest are

$$\gamma_1 = \frac{1}{\lambda_1 \tau_1}, \qquad \gamma_2 = \frac{1}{\lambda_2 \tau_2} \quad \text{and} \quad \gamma = \frac{1}{\lambda_1 \tau_1 + \lambda_2 \tau_2}.$$

First, we check assumption (15). Using the above approximation of $H^{\leftarrow}$, it follows, for $s \leq t < 1$ and $s$ large enough, that

$$|p(H^{\leftarrow}(t)) - p(H^{\leftarrow}(s))| \leq \widetilde{C}((1-s)^{\gamma\tau} - (1-t)^{\gamma\tau})$$

for some $\widetilde{C} > 0$. It now easily follows that in the case $\gamma\tau \geq 1$, the left-hand side of (15) tends to 0. In the case $\gamma\tau < 1$, the left-hand side of (15) is of order $\sqrt{k}(\frac{\sqrt{k}}{n})^{\gamma\tau} = \sqrt{k}(\frac{n}{k})^{\rho}k^{\rho/2}$, which tends to 0 when (14) holds (see below).

The asymptotic bias of $\sqrt{k}(\widehat{\gamma}_{Z,k,n}^{(\cdot)} - \gamma)$ can be explicitly computed (from Corollary 1) and is asymptotically equivalent to

$$-\eta(\beta_1^{\lambda_1}\beta_2^{\lambda_2})^{\rho}\sqrt{k}\left(\frac{n}{k}\right)^{\rho} \cdot \begin{cases} \dfrac{\gamma\rho}{1-\rho}, & \text{for the Hill estimator,} \\[2mm] \dfrac{\rho(1+\gamma)(\gamma+\rho)}{(1-\rho)(1-\rho+\gamma)}, & \text{for the } ML\text{-estimator,} \\[2mm] \dfrac{\rho[\rho+\gamma(1-\rho)]}{(1-\rho)^2}, & \text{for the moment and } UH\text{-estimators.} \end{cases}$$

They are all of the same order.

We obtain another bias term from assumption (14). Direct computations, using (12) and $p = \frac{\gamma_2}{\gamma_1 + \gamma_2}$, lead to

$$p(z) - p = \frac{\gamma^2}{\gamma_1 \gamma_2}[-\beta_1 z^{-\tau_1}(1 + o(1)) + \beta_2 z^{-\tau_2}(1 + o(1))]$$

when $\tau_1 \neq \tau_2$, or $\tau_1 = \tau_2$ and $\beta_1 \neq \beta_2$. Consequently, assumption (14) is equivalent to

$$\beta \frac{\gamma^2}{\gamma_1 \gamma_2}(\beta_1^{\lambda_1}\beta_2^{\lambda_2})^{\rho}\frac{1}{1-\rho}\sqrt{k}\left(\frac{n}{k}\right)^{\rho} \longrightarrow \alpha_2,$$

with

$$\beta = \begin{cases} -\beta_1, & \text{if } \tau_1 < \tau_2, \\ \beta_2, & \text{if } \tau_1 > \tau_2, \\ \beta_2 - \beta_1, & \text{if } \tau_1 = \tau_2. \end{cases}$$

So both bias terms are of the same order. Only when $\tau_1 = \tau_2$ and $\beta_1 = \beta_2$ (in particular, when $F \equiv G$) the biases of the estimators of $\gamma$ dominate.

**Example 2.** $X \sim$ reverse Burr$(\beta_1, \tau_1, \lambda_1, x_+)$ and $Y \sim$ reverse Burr$(\beta_2, \tau_2, \lambda_2, x_+)$, $\beta_1, \tau_1, \lambda_1, \beta_2, \tau_2, \lambda_2, x_+ > 0$.

In that case,

$$1 - F(x) = \left(\frac{\beta_1}{\beta_1 + (x_+ - x)^{-\tau_1}}\right)^{\lambda_1}$$

$$= (x_+ - x)^{\tau_1 \lambda_1} \beta_1^{\lambda_1} (1 + \beta_1 (x_+ - x)^{\tau_1})^{-\lambda_1}, \qquad x < x_+;$$

$$1 - G(x) = \left(\frac{\beta_2}{\beta_2 + (x_+ - x)^{-\tau_2}}\right)^{\lambda_2}$$

$$= (x_+ - x)^{\tau_2 \lambda_2} \beta_2^{\lambda_2} (1 + \beta_2 (x_+ - x)^{\tau_2})^{-\lambda_2}, \qquad x < x_+.$$

Define $\tau$ and $\eta$ as in Example 1, but now set $\rho = \gamma \tau$. We can infer that

$$U_H(x) = H^{\leftarrow}\left(1 - \frac{1}{x}\right) = x_+ - (\beta_1^{\lambda_1} \beta_2^{\lambda_2} x)^{-1/(\tau_1 \lambda_1 + \tau_2 \lambda_2)} [1 - \gamma \eta (\beta_1^{\lambda_1} \beta_2^{\lambda_2} x)^{\rho} (1 + o(1))].$$

The parameters of interest are

$$\gamma_1 = -\frac{1}{\lambda_1 \tau_1}, \qquad \gamma_2 = -\frac{1}{\lambda_2 \tau_2}, \qquad \gamma = -\frac{1}{\lambda_1 \tau_1 + \lambda_2 \tau_2} \quad \text{and} \quad \tau_F = \tau_G = \tau_H = x_+.$$

Note that we can easily prove (as in Example 1) that assumption (15) is satisfied if we assume (14).

The asymptotic bias of $\sqrt{k}(\widehat{\gamma}_{Z,k,n}^{(\cdot)} - \gamma)$ can be explicitly computed (again, from Corollary 1) and is asymptotically equivalent to

- for the *UH*-estimator:

$$\begin{cases} -\dfrac{\gamma^2 \tau (1 - \gamma)(1 + \tau)}{(1 - \gamma - \gamma \tau)(1 - \gamma \tau)} \eta (\beta_1^{\lambda_1} \beta_2^{\lambda_2})^{\rho} \sqrt{k} \left(\dfrac{n}{k}\right)^{\rho}, & \text{if } \tau < 1, \\[4mm] \dfrac{\gamma^2}{(1 - \gamma)(1 - 2\gamma)} (\beta_1^{\lambda_1} \beta_2^{\lambda_2})^{\rho} \left[-2\eta(1 - \gamma) + \dfrac{1}{x_+}\right] \sqrt{k} \left(\dfrac{n}{k}\right)^{\rho}, & \text{if } \tau = 1, \\[4mm] \dfrac{\gamma^2}{(1 - \gamma)(1 - 2\gamma) x_+} (\beta_1^{\lambda_1} \beta_2^{\lambda_2})^{\gamma} \sqrt{k} \left(\dfrac{n}{k}\right)^{\gamma}, & \text{if } \tau > 1; \end{cases}$$

- for the moment estimator:

$$\begin{cases} -\dfrac{\gamma^2 \tau (1 - \gamma)(1 + \tau)(1 - 2\gamma)}{(1 - \gamma - \gamma \tau)(1 - 2\gamma - \gamma \tau)} \eta (\beta_1^{\lambda_1} \beta_2^{\lambda_2})^{\rho} \sqrt{k} \left(\dfrac{n}{k}\right)^{\rho}, & \text{if } \tau < 1, \\[4mm] -\dfrac{\gamma^2}{(1 - \gamma)(1 - 3\gamma)} (\beta_1^{\lambda_1} \beta_2^{\lambda_2})^{\rho} \left[2\eta(1 - \gamma)^2 - \dfrac{\gamma + 1}{x_+}\right] \sqrt{k} \left(\dfrac{n}{k}\right)^{\rho}, & \text{if } \tau = 1, \\[4mm] -\dfrac{\gamma}{(1 - \gamma) x_+} (\beta_1^{\lambda_1} \beta_2^{\lambda_2})^{\gamma} \sqrt{k} \left(\dfrac{n}{k}\right)^{\gamma}, & \text{if } \tau > 1; \end{cases}$$

- for the *ML*-estimator, if $\gamma > -\frac{1}{2}$:

$$-\frac{\gamma^2 \tau (1 + \gamma)(1 + \tau)}{(1 - \gamma \tau)(1 + \gamma - \gamma \tau)} \eta (\beta_1^{\lambda_1} \beta_2^{\lambda_2})^{\rho} \sqrt{k} \left(\frac{n}{k}\right)^{\rho}.$$

They are all of the same order if $\tau \leq 1$, otherwise the biases of the moment and *UH*-estimators dominate that of the *ML*-estimator.

Similarly to Example 1, if $\tau_1 \neq \tau_2$, or $\tau_1 = \tau_2$ and $\beta_1 \neq \beta_2$, direct computations lead to

$$p(z) - p = \frac{\gamma^2}{\gamma_1 \gamma_2}[-\beta_1(x_+ - z)^{\tau_1}(1 + o(1)) + \beta_2(x_+ - z)^{\tau_2}(1 + o(1))].$$

Consequently, assumption (14) is equivalent, in that case, to

$$\beta \frac{\gamma^2}{\gamma_1 \gamma_2}(\beta_1^{\lambda_1} \beta_2^{\lambda_2})^\rho \frac{1}{1 - \rho}\sqrt{k}\left(\frac{n}{k}\right)^\rho \longrightarrow \alpha_2.$$

Again, this order is the same as the order of the asymptotic bias terms of all of the estimators in case $\tau \leq 1$ and dominated by the one of the moment and *UH*-estimators otherwise. When $\tau_1 = \tau_2$ and $\beta_1 = \beta_2$, the biases of the estimators of $\gamma$ dominate.

**Example 3.** *$X, Y \sim$ logistic.*

In that case,

$$1 - F(x) = 1 - G(x) = \frac{2}{1 + e^x}, \qquad x > 0.$$

Hence,

$$U_H(x) = \log(2\sqrt{x} - 1).$$

We have $\gamma_1 = \gamma_2 = \gamma = 0$. Since $F \equiv G$, we immediately obtain $p(\cdot) \equiv \frac{1}{2}$ and $\alpha_2 = 0$.

According to Corollary 1, the asymptotic bias of $\sqrt{k}(\widehat{\gamma}_{Z,k,n}^{(\cdot)} - 0)$ is asymptotically equivalent to $\frac{\sqrt{k}}{\log^2 n/k}$ for the *UH*- or the moment estimator and to $-\frac{1}{9}\frac{k}{\sqrt{n}}$ for the *ML*-estimator.

In order to illustrate these three examples, we simulate 100 samples of size 500 from the following distributions:

- a Burr$(10, 4, 1)$ censored by a Burr$(10, 1, 0.5)$;
- a reverse Burr$(1, 8, 0.5, 10)$ censored by a reverse Burr$(10, 1, 0.5, 10)$;
- a logistic censored by a logistic.

For the first two examples, $p = \frac{8}{9}$, meaning that the percentage of censoring in the right tail is close to 11%. In the last example, $p(\cdot) \equiv p = \frac{1}{2}$, that is, the percentage of censoring is as high as 50%. In the first case, we have $\gamma_1 = \frac{1}{4}, \gamma = \frac{2}{9}$ and $\rho = -\frac{2}{9}$, in the second case $\gamma_1 = -\frac{1}{4}, \gamma = -\frac{2}{9}$ and, again, $\rho = -\frac{2}{9}$. In the third example, $\gamma_1 = \gamma = 0$. In all three examples, panels (a) and (c) (in Figures 2–4) represent the median for the index and the extreme quantile, respectively, whereas panels (b) and (d) represent the empirical mean square errors (MSE) based on the 100 samples. The small value of $\varepsilon$ is $\frac{1}{50}$. All of these plotted estimators are adapted to censoring. The horizontal line represents the true value of the parameter.

In the first example, we can observe, in the case of the estimation of the index, the superiority of the Hill estimator adapted to censoring in terms of MSE, the three others being quite similar. For the extreme quantile estimators, however, there is much less to decide between all of the estimators: they are very stable and close to the true value of the parameter. A similar observation can be made for the second and third examples, with a slight advantage for the *UH*-estimator, only in the case of the estimation of the index.

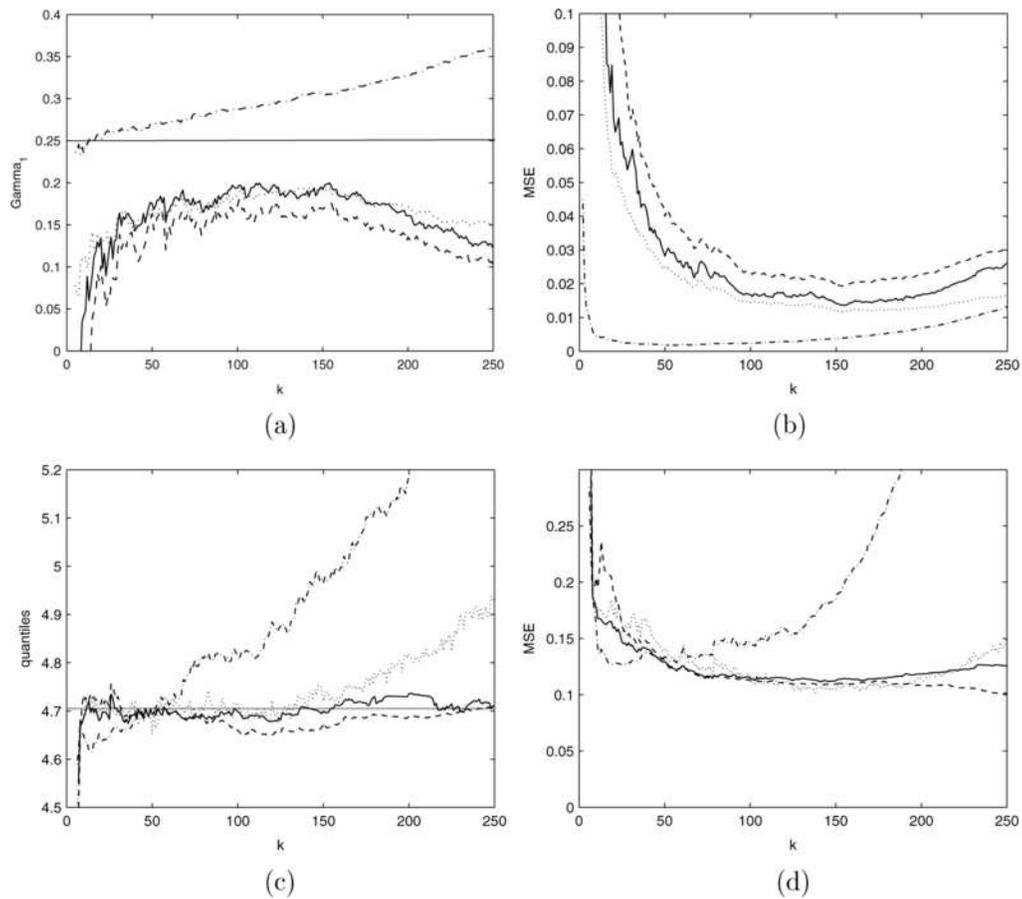

**Figure 2.** A Burr$(10, 4, 1)$ distribution censored by a Burr$(10, 1, 0.5)$ distribution: *UH*-estimator (dotted line), moment estimator (full line), *ML*-estimator (dashed line) and Hill estimator (dashed-dotted line); (a) median and (b) MSE for the extreme value index; (c) median and (d) MSE for the extreme quantile with $\varepsilon = \frac{1}{50}$.

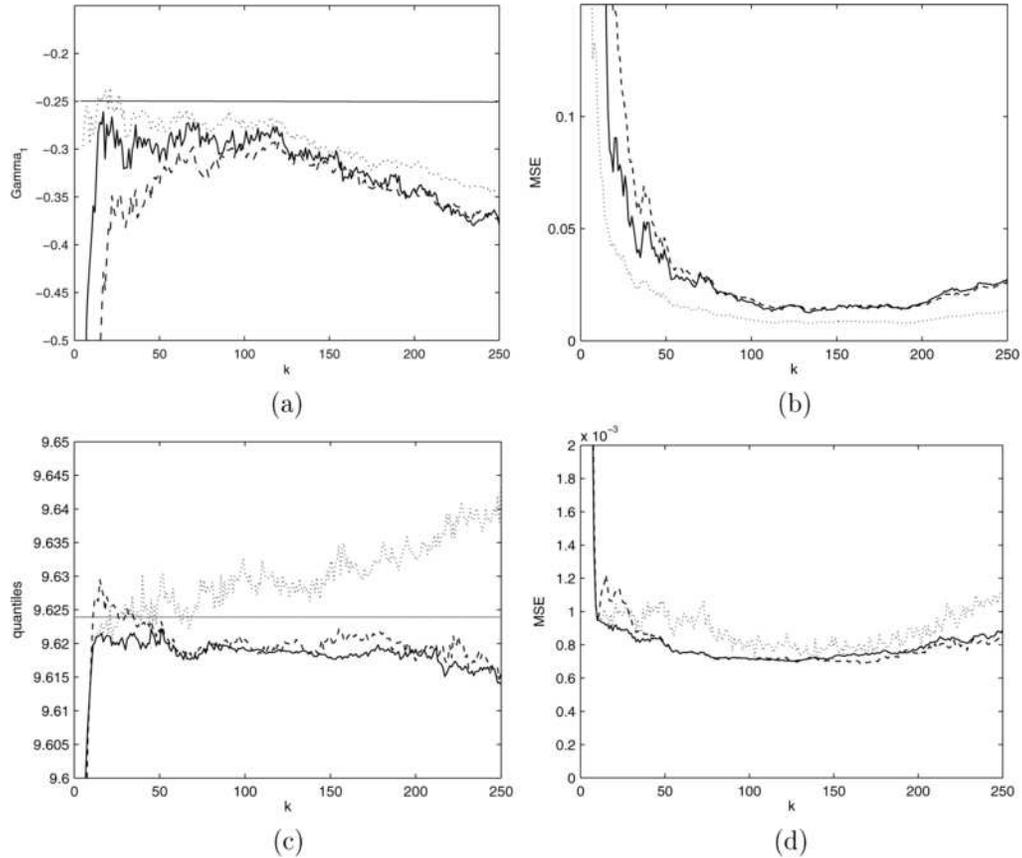

**Figure 3.** A reverse Burr$(1, 8, 0.5, 10)$ distribution censored by a reverse Burr$(10, 1, 0.5, 10)$ distribution: *UH*-estimator (dotted line), moment estimator (full line), *ML*-estimator (dashed line); (a) median and (b) MSE for the extreme value index; (c) median and (d) MSE for the extreme quantile with $\varepsilon = \frac{1}{50}$.

## 4. Application to AIDS survival data

We return to our real data set presented in Section 1 and used in Section 2, that is, the Australian AIDS survival data for the male patients diagnosed before 1 July 1991. The sample size is 2754.

First, we estimate $p = \lim_{z \to \tau_H} p(z)$. In Figure 5, we see $\widehat{p}$ as a function of $k$. Clearly, there is a stable part in the plot when $k$ ranges from about 75 to 175; for higher $k$, the bias sets in. Note that $\widehat{p}$ is the mean of 0–1 variables, so for a sample of this size, the estimator is already very accurate. Therefore, we estimate $p$ with the corresponding vertical level in the plot, which is 0.28.

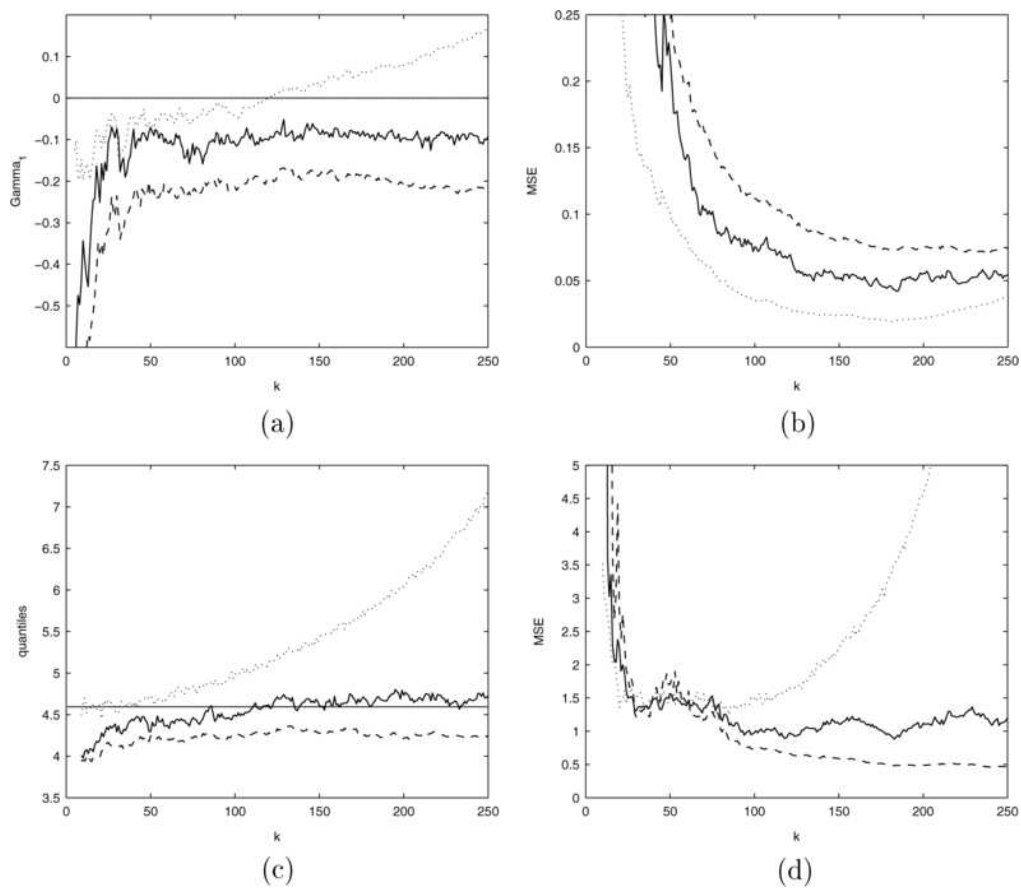

**Figure 4.** A logistic distribution censored by a logistic distribution: *UH*-estimator (dotted line), moment estimator (full line), *ML*-estimator (dashed line); (a) median and (b) MSE for the extreme value index; (c) median and (d) MSE for the extreme quantile with $\varepsilon = \frac{1}{50}$.

We now continue with the estimation of the extreme value index $\gamma_1$ and an extreme quantile $F^{\leftarrow}(1-\varepsilon)$, using the *UH*-method (as in Section 2). We will again plot these estimators as functions of $k$, but already replacing $\widehat{p} = \widehat{p}(k)$ with its estimate 0.28 in order to prevent that the bias plays a dominant role for values of $k$ larger than 200, say.

In Figure 6(a), the estimator of the extreme value index is presented, whereas Figure 6(b) shows the extreme quantile estimator for $\varepsilon = 0.001$. The estimator of $\gamma_1$ is quite stable for values of $k$ between 200 and 300. We estimate it with 0.14. We estimate the extreme quantile with $k$ values in the same range because that range again gives a stable part in the plot. The corresponding estimated survival time is as high as about 25 years. So, although the estimated median survival time has the low value 1.3 years, we find that exceptionally strong males can survive AIDS for 25 years.

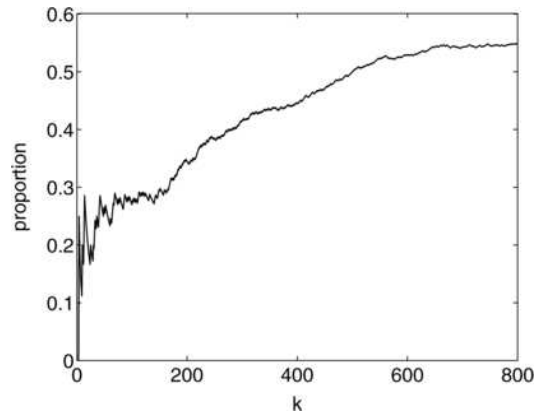

**Figure 5.** Estimator of $p$ for the Australian AIDS survival data for the male patients.

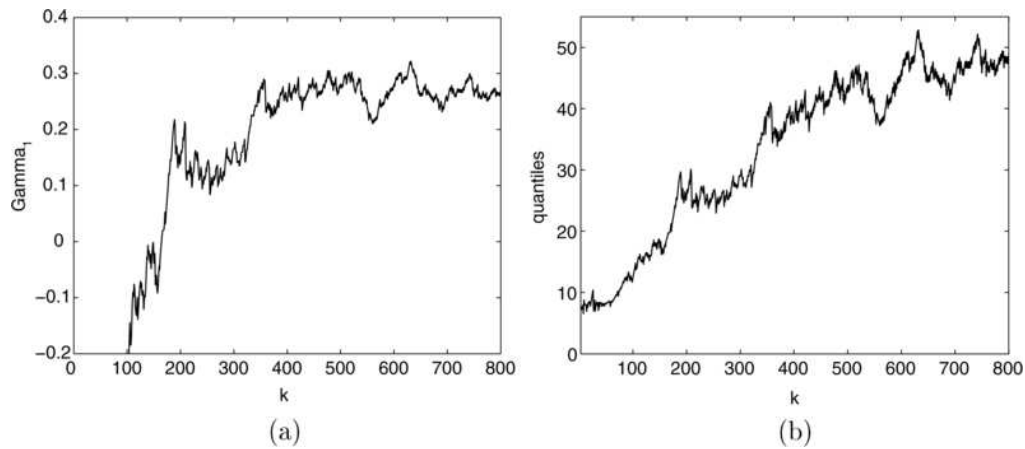

**Figure 6.** *UH*-estimator (a) for the extreme value index and (b) for the extreme quantile with $\varepsilon = 0.001$, for the Australian AIDS survival data for the male patients.

# Acknowledgements

We are grateful to two referees for their thoughtful and constructive comments.